\documentclass[reqno]{amsart}

\newtheorem{Theorem}{Theorem}[section]
\newtheorem{Proposition}[Theorem]{Proposition}

\theoremstyle{remark}

\numberwithin{equation}{section}

\allowdisplaybreaks

\begin{document}

\title{Some more semi-finite forms of bilateral basic hypergeometric series}

\author[Fr\'ed\'eric Jouhet]{Fr\'ed\'eric Jouhet$^*$}
\address{Institut Girard Desargues, Universit\'e Claude Bernard (Lyon 1),
69622 Villeurbanne Cedex, France}
\email{jouhet@igd.univ-lyon1.fr}
\urladdr{http://igd.univ-lyon1.fr/home/jouhet}

\thanks{$^*$Partially supported by EC's IHRP Programme,
grant HPRN-CT-2001-00272, ``Algebraic Combinatorics in Europe''.}


\begin{abstract}
We prove some new semi-finite forms of bilateral basic hypergeometric series. One of them yields in a direct limit Bailey's celebrated ${}_6\psi_6$ summation formula, answering a question recently raised by Chen and Fu ({\em Semi-Finite Forms of Bilateral Basic Hypergeometric Series}, Proc. Amer. Math. Soc., to appear).
\end{abstract}

\maketitle

\section{Introduction}
There is a standard method for obtaining a bilateral identity from a \emph{unilateral terminating identity}, which was already utilized by Cauchy~\cite{Ca} in his second proof of Jacobi's~\cite{Jac} famous triple product identity. The idea of this method is to start from a finite unilateral summation and to shift the index of summation, say $k$ ($0\leq k\leq 2n$), by $n$ :
\begin{equation}\label{finite}
\sum_{k=0}^{2n}a(k)=\sum_{k=-n}^{n}a(k+n),
\end{equation}
and then let $n\to\infty$ whenever it is possible after some manipulations. The same method has also been exploited by Bailey~\cite[Secs.~3 and 6]{Ba1}, \cite{Ba2}, Slater~\cite[Sec.~6.2]{Sl}, Schlosser~\cite{Sc} and Schlosser and the author~\cite{JS}.\\ 
Recently, Chen and Fu \cite{CF} used a method different from the previous one, as they started from \emph{unilateral infinite summations} to derive \emph{semi-finite forms of bilateral basic hypergeometric series}. The process can be summarized as follows~:
\begin{equation}\label{infinite}
\sum_{k\geq 0}a(k)=\sum_{k\geq -n}a(k+n),
\end{equation}
and then let $n\to\infty$ whenever it is possible after some manipulations. The right-hand side of (\ref{finite}) (resp. (\ref{infinite})) can be seen as a finite (resp. semi-finite) form of a bilateral series. Chen and Fu have found in \cite{CF} semi-finite forms of Ramanujan's ${}_1\psi_1$ summation formula (cf. \cite[Appendix (II.29)]{GR}), of a 
${}_2\psi_2$ formula due to Bailey \cite[Ex. 5.20(i)]{GR}, and of Bailey's~\cite[Eq.~(4.7)]{Ba1} ${}_6\psi_6$ summation formula (\cite[Appendix (II.33)]{GR}), which can be written as follows~: 
\begin{multline}\label{6psi6}
{}_6\psi_6\!\left[\begin{matrix}q\sqrt{a},-q\sqrt{a},b,c,d,e\\
\sqrt{a},-\sqrt{a},aq/b,aq/c,aq/d,aq/e\end{matrix};
q,\frac{qa^2}{bcde}\right]\\
=\frac{(q,aq,q/a,aq/bc,aq/bd,aq/be,aq/cd,aq/ce,aq/de)_\infty}
{(q/b,q/c,q/d,q/e,aq/b,aq/c,aq/d,aq/e,qa^2/bcde)_\infty}
\end{multline}
(see the end of this introduction for the notations), where $|q|<1$ and $|qa^2/bcde|<1$.\\
At the end of \cite{CF}, Chen and Fu mention the problem of finding a proof of (\ref{6psi6}) using a semi-finite (or even finite) form which would yield (\ref{6psi6}) in a direct limit. Indeed, after letting $n\to\infty$ in their semi-finite form of (\ref{6psi6}), one needs to use Ramanujan's ${}_1\psi_1$ summation formula to derive (\ref{6psi6}).\\
In this paper, we use the method developped in \cite{CF} to find, among other results, a new semi-finite form of (\ref{6psi6}) which answers the question raised by Chen and Fu.\\
After explaining some notations in the end of this introduction, we show in section~2 how the method in \cite{CF} can be applied to yield in a direct limit (\ref{6psi6}), starting from a nonterminating extension of Jackson's formula due to Bailey~\cite[Appendix (II.25)]{GR}. We give two other applications of this method in section~3, which yield in a direct limit a transformation formula for a $_6\psi_6$ series proved in \cite{JS} and a transformation formula for a $_8\psi_8$ series in terms of two $_8\phi_7$ series and a $_8\psi_8$ series, which seems to be new.

Other proofs of Bailey's very-well-poised ${}_6\psi_6$ summation
had been given by Bailey~\cite{Ba1}, Slater and Lakin~\cite{SL},
Andrews~\cite{And}, Askey and Ismail~\cite{AI},
Askey~\cite{AII}, Schlosser~\cite{Sc6psi6} and Schlosser and the author~\cite{JS}. It is worth noting that the elegant proof of Askey and Ismail in~\cite{AI} uses an argument of \emph{analytic continuation} together with the shift (\ref{infinite}), but used from right to left.

{\bf Notation:}
It is appropriate to recall some standard notations for \emph{$q$-series} and
\emph{basic hypergeometric series}.

Let $q$ be a fixed complex parameter (the ``base'') with $0<|q|<1$.
The \emph{$q$-shifted factorial} is defined for any complex
parameter $a$ by
\begin{equation*}
(a)_\infty\equiv (a;q)_\infty:=\prod_{j\geq 0}(1-aq^j)
\end{equation*}
and
\begin{equation*}
(a)_k\equiv (a;q)_k:=\frac{(a;q)_\infty}{(aq^k;q)_\infty},
\end{equation*}
where $k$ is any integer.
Since the same base $q$ is used throughout this paper,
it may be readily omitted (in notation, writing $(a)_k$ instead of $(a;q)_k$, etc) which will not lead
to any confusion. For brevity, write
\begin{equation*}
(a_1,\ldots,a_m)_k:=(a_1)_k\cdots(a_m)_k,
\end{equation*}
where $k$ is an integer or infinity.
Further, recall the definition of \emph{basic hypergeometric series},
\begin{equation*}
{}_s\phi_{s-1}\!\left[\begin{matrix}a_1,\dots,a_s\\
b_1,\dots,b_{s-1}\end{matrix};q,z\right]:=
\sum_{k=0}^\infty\frac{(a_1,\dots,a_s)_k}{(q,b_1,\dots,b_{s-1})_k}z^k,
\end{equation*}
and of \emph{bilateral basic hypergeometric series},
\begin{equation*}
{}_s\psi_s\!\left[\begin{matrix}a_1,\dots,a_s\\
b_1,\dots,b_s\end{matrix};q,z\right]:=
\sum_{k=-\infty}^\infty\frac{(a_1,\dots,a_s)_k}{(b_1,\dots,b_s)_k}z^k.
\end{equation*}

See Gasper and Rahman's text~\cite{GR} for a
comprehensive study of the theory of basic hypergeometric series.
In particular, the computations in this paper rely on some elementary
identities for $q$-shifted factorials, listed in \cite[Appendix~I]{GR}.

\section{A new semi-finite form of Bailey's ${}_6\psi_6$ summation formula}
Consider Bailey's nonterminating extension of Jackson's ${}_{8}\phi_7$ summation~\cite[Appendix (II.25)]{GR}
\begin{multline}\label{8phi7ext}
{}_{8}\phi_7\!\left[\begin{matrix}a,\,q\sqrt{a},-q\sqrt{a},b,c,d,e,f\\
\sqrt{a},-\sqrt{a},aq/b,aq/c,aq/d,aq/e,aq/f\end{matrix};q,q\right]\\
=\frac{b}{a}\frac{(aq,c,d,e,f,bq/a,bq/c,bq/d,bq/e,bq/f)_\infty}{(aq/b,aq/c,aq/d,aq/e,aq/f,bc/a,bd/a,be/a,bf/a,b^2q/a)_\infty}\\\times
{}_{8}\phi_7\!\left[\begin{matrix}b^2/a,\,
qb/\sqrt{a},-qb/\sqrt{a},b,bc/a,bd/a,be/a,bf/a\\
b/\sqrt{a},-b/\sqrt{a},bq/a,bq/c,bq/d,bq/e,bq/f\end{matrix};q,q\right]\\
+\frac{(aq,b/a,aq/cd,aq/ce,aq/cf,aq/de,aq/df,aq/ef)_\infty}{(aq/c,aq/d,aq/e,aq/f,bc/a,bd/a,be/a,bf/a)_\infty},
\end{multline}
where $qa^2=bcdef$.\\
 Note that (\ref{8phi7ext}) can be proved by specializing $qa^2=bcdef$ in Bailey's 3-term transformation formula for a nonterminating very-well-poised ${}_{8}\phi_7$~\cite[Appendix (III.37)]{GR}, which was the starting point in \cite{CF} to prove (\ref{6psi6}), and then using the sum of a very-well-poised ${}_{6}\phi_5$~\cite[Appendix (II.20)]{GR}.\\
 Now, using (\ref{8phi7ext}), we can derive the following semi-finite form of (\ref{6psi6}).
\begin{Proposition}
\begin{multline}\label{semifinite6psi6}
\sum_{k\geq -n}\frac{(aq^{-n},q\sqrt{a},-q\sqrt{a},bq^n,c,d,e,f)_k}{(q^{1+n},\sqrt{a},-\sqrt{a},aq^{1-n}/b,aq/c,aq/d,aq/e,aq/f)_k}q^k\\
=\frac{(aq,c,d,e,f,bq^{1+2n}/a,bq^{1+n}/c,bq^{1+n}/d,bq^{1+n}/e,bq^{1+n}/f)_\infty}{(aq/b,aq/c,aq/d,aq/e,aq/f,bcq^n/a,bdq^n/a,beq^n/a,bfq^n/a,b^2q^{1+2n}/a)_\infty}\\
\times\frac{b^{n+1}}{a}\frac{(q,q/a)_n}{(b,b/a)_n}\\
\times{}_{8}\phi_7\!\left[\begin{matrix}b^2q^{2n}/a,\,
bq^{1+n}/\sqrt{a},-bq^{1+n}/\sqrt{a},b,bcq^n/a,bdq^n/a,beq^n/a,bfq^n/a\\
bq^n/\sqrt{a},-bq^n/\sqrt{a},bq^{1+2n}/a,bq^{1+n}/c,bq^{1+n}/d,bq^{1+n}/e,bq^{1+n}/f\end{matrix};q,q\right]\\
+\frac{(aq,aq/cd,aq/ce,aq/cf,aq/de,aq/df,aq/ef,bq^{n}/a)_\infty}{(aq/c,aq/d,aq/e,aq/f,bcq^n/a,bdq^n/a,beq^n/a,bfq^n/a)_\infty}\\
\times\frac{(q,q/a)_n}{(b,q/c,q/d,q/e,q/f)_n},
\end{multline}
where $b=qa^2/cdef$.
\end{Proposition}
\begin{proof}
By shifting the index of summation by $n$, the left-hand side of (\ref{8phi7ext}) is equal to
\begin{multline*}
\frac{1-aq^{2n}}{1-a}\frac{(a,b,c,d,e,f)_n}{(q,aq/b,aq/c,aq/d,aq/e,aq/f)_n}q^n\\
\times\sum_{k\geq -n}\frac{(aq^{n},\sqrt{a}q^{1+n},-\sqrt{a}q^{1+n},bq^n,cq^n,dq^n,eq^n,fq^n)_k}{(q^{1+n},\sqrt{a}q^n,-\sqrt{a}q^n,aq/b,aq^{1+n}/c,aq^{1+n}/d,aq^{1+n}/e,aq^{1+n}/f)_k}q^k.
\end{multline*}
Next, on both sides of (\ref{8phi7ext}), replace $a$, $c$, $d$, $e$ and $f$ by $aq^{-2n}$, $cq^{-n}$, $dq^{-n}$, $eq^{-n}$ and $fq^{-n}$ respectively. Note that the condition $qa^2=bcdef$ is equivalent to $b=qa^2/cdef$, thus $b$ remains unchanged. We get
\begin{multline*}
\frac{1-a}{1-aq^{-2n}}\frac{(aq^{-2n},b,cq^{-n},dq^{-n},eq^{-n},fq^{-n})_n}{(q,aq^{1-2n}/b,aq^{1-n}/c,aq^{1-n}/d,aq^{1-n}/e,aq^{1-n}/f)_n}q^n\\
\times\sum_{k\geq -n}\frac{(aq^{-n},q\sqrt{a},-q\sqrt{a},bq^n,c,d,e,f)_k}{(q^{1+n},\sqrt{a},-\sqrt{a},aq^{1-n}/b,aq/c,aq/d,aq/e,aq/f)_k}q^k\\
=\frac{b}{aq^{-2n}}\frac{(aq^{1-2n},cq^{-n},dq^{-n},eq^{-n},fq^{-n})_\infty}{(aq^{1-2n}/b,aq^{1-n}/c,aq^{1-n}/d,aq^{1-n}/e,aq^{1-n}/f)_\infty}\\
\times\frac{(bq^{1+2n}/a,bq^{1+n}/c,bq^{1+n}/d,bq^{1+n}/e,bq^{1+n}/f)_\infty}{(bcq^n/a,bdq^n/a,beq^n/a,bfq^n/a,b^2q^{1+2n}/a)_\infty}\\
\times{}_{8}\phi_7\!\left[\begin{matrix}b^2q^{2n}/a,\,
bq^{1+n}/\sqrt{a},-bq^{1+n}/\sqrt{a},b,bcq^n/a,bdq^n/a,beq^n/a,bfq^n/a\\
bq^n/\sqrt{a},-bq^n/\sqrt{a},bq^{1+2n}/a,bq^{1+n}/c,bq^{1+n}/d,bq^{1+n}/e,bq^{1+n}/f\end{matrix};q,q\right]\\
+\frac{(aq^{1-2n},bq^{2n}/a,aq/cd,aq/ce,aq/cf,aq/de,aq/df,aq/ef)_\infty}{(aq^{1-n}/c,aq^{1-n}/d,aq^{1-n}/e,aq^{1-n}/f,bcq^n/a,bdq^n/a,beq^n/a,bfq^n/a)_\infty}.
\end{multline*}
This can be rewritten as
\begin{multline}\label{c}
\sum_{k\geq -n}\frac{(aq^{-n},q\sqrt{a},-q\sqrt{a},bq^n,c,d,e,f)_k}{(q^{1+n},\sqrt{a},-\sqrt{a},aq^{1-n}/b,aq/c,aq/d,aq/e,aq/f)_k}q^k\\
=\frac{1-aq^{-2n}}{1-a}\frac{bq^n}{a}\frac{(q,aq^{1-2n}/b)_n}{(b,aq^{-2n})_n}\frac{(aq^{1-2n})_\infty}{(aq^{1-2n}/b)_\infty}\\
\times\frac{(c,d,e,f,bq^{1+2n}/a,bq^{1+n}/c,bq^{1+n}/d,bq^{1+n}/e,bq^{1+n}/f)_\infty}{(aq/c,aq/d,aq/e,aq/f,bcq^n/a,bdq^n/a,beq^n/a,bfq^n/a,b^2q^{1+2n}/a)_\infty}\\
\times{}_{8}\phi_7\!\left[\begin{matrix}b^2q^{2n}/a,\,
bq^{1+n}/\sqrt{a},-bq^{1+n}/\sqrt{a},b,bcq^n/a,bdq^n/a,beq^n/a,bfq^n/a\\
bq^n/\sqrt{a},-bq^n/\sqrt{a},bq^{1+2n}/a,bq^{1+n}/c,bq^{1+n}/d,bq^{1+n}/e,bq^{1+n}/f\end{matrix};q,q\right]\\
+\frac{1-aq^{-2n}}{1-a}q^{-n}\frac{(q,aq^{1-2n}/b)_n(aq^{1-2n})_\infty}{(aq^{-2n},b,cq^{-n},dq^{-n},eq^{-n},fq^{-n})_n}\\
\times\frac{(bq^{2n}/a,aq/cd,aq/ce,aq/cf,aq/de,aq/df,aq/ef)_\infty}{(aq/c,aq/d,aq/e,aq/f,bcq^n/a,bdq^n/a,beq^n/a,bfq^n/a)_\infty}.
\end{multline}
Now we use the three following elementary identities to simplify the right-hand side of (\ref{c})~:
\begin{eqnarray}
\frac{(xq^{-2n})_\infty}{(xq^{-2n})_n}&=&(-1)^nx^nq^{-(n^2+n)/2}(q/x)_n(x)_\infty,\label{elementary1}\\
(xq^{-2n})_n&=&(-1)^nx^nq^{-(3n^2+n)/2}(q^{n+1}/x)_n\label{elementary2},\\
(xq^{-n})_n&=&(-1)^nx^nq^{-(n^2+n)/2}(q/x)_n,\label{elementary3}
\end{eqnarray}
and we obtain (\ref{semifinite6psi6}) after simplifications.
\end{proof}
Now, one may let $n\to\infty$ in (\ref{semifinite6psi6}), assuming $|qa^2/cdef|<1$ (i.e. $|b|<1$), while appealing to Tannery's theorem~\cite{Br} for being able to interchange limit and summation. As the first term on the right-hand side of (\ref{semifinite6psi6}) tends to 0, this gives immediately Bailey's ${}_6\psi_6$ summation formula (\ref{6psi6}) with $b$ replaced by $f$.

\section{Other consequences}
We give in this section two other application of the previous process. Consider first the following transformation formula for a non terminating very-well-poised ${}_8\phi_7$ series~\cite[Appendix (III.23)]{GR} 
\begin{multline}\label{8phi7}
{}_{8}\phi_7\!\left[\begin{matrix}a,\,q\sqrt{a},-q\sqrt{a},b,c,d,e,f\\
\sqrt{a},-\sqrt{a},aq/b,aq/c,aq/d,aq/e,aq/f\end{matrix};q,\frac{q^2a^2}{bcdef}\right]\\
=\frac{(aq,aq/ef,\lambda q/e,\lambda q/f)_\infty}
{(aq/e,aq/f,\lambda q/ef,\lambda q)_\infty}\\\times
{}_{8}\phi_7\!\left[\begin{matrix}\lambda,\,
q\sqrt{\lambda},-q\sqrt{\lambda},\lambda
b/a,\lambda c/a,\lambda d/a,e,f\\
\sqrt{\lambda},-\sqrt{\lambda},aq/b,aq/c,aq/d,\lambda q/e,
\lambda q/f\end{matrix};q,\frac{aq}{ef}\right],
\end{multline}
where $\lambda=qa^2/bcd$, $|q^2a^2/bcdef|<1$ and $|aq/ef|<1$.\\
 Note that (\ref{8phi7}) is nothing else but the $n\to\infty$ case of Bailey's~\cite{Ba1} transformation formula for a very-well-poised ${}_{10}\phi_9$ series~\cite[Appendix (III.28)]{GR}, which was the starting point in \cite{JS} for the derivation of (\ref{6psi6}). Now, using (\ref{8phi7}), we can prove the following semi-finite identity.
\begin{Proposition}
\begin{multline}\label{semifinite}
\sum_{k\geq -n}\frac{(aq^{-n},q\sqrt{a},-q\sqrt{a},bq^n,c,d,e,f)_k}{(q^{1+n},\sqrt{a},-\sqrt{a},aq^{1-n}/b,aq/c,aq/d,aq/e,aq/f)_k}\left(\frac{q^2a^2}{bcdef}\right)^k\\
=\frac{(aq,aq/ef,\lambda q/e,\lambda q/f)_\infty}{(aq/e,aq/f,\lambda q/ef,\lambda q)_\infty}\frac{(\lambda b/a,q/a,aq/\lambda c,aq/\lambda d)_n}{(b,q/\lambda,q/c,q/d)_n}\\
\times\sum_{k\geq -n}\frac{(\lambda q^{-n},q\sqrt{\lambda},-q\sqrt{\lambda},\lambda bq^{n}/a,\lambda c/a,\lambda d/a,e,f)_k}{(q^{1+n},\sqrt{\lambda},-\sqrt{\lambda},aq^{1-n}/b,aq/c,aq/d,\lambda q/e,\lambda q/f)_k}\left(\frac{aq}{ef}\right)^k,
\end{multline}
where $\lambda=qa^2/bcd$ and $|q^2a^2/bcdef|<1$.
\end{Proposition}
\begin{proof}
By shifting the index of summation by $n$ on both sides of (\ref{8phi7}), we get 
\begin{multline}\label{a}
\frac{1-aq^{2n}}{1-a}\frac{(a,b,c,d,e,f)_n}{(q,aq/b,aq/c,aq/d,aq/e,aq/f)_n}\left(\frac{q^2a^2}{bcdef}\right)^n\\
\times\sum_{k\geq -n}\frac{(aq^n,\sqrt{a}q^{1+n},-\sqrt{a}q^{1+n},bq^n,cq^n,dq^n)_k}{(q^{1+n},\sqrt{a}q^{n},-\sqrt{a}q^{n},aq^{1+n}/b,aq^{1+n}/c,aq^{1+n}/d)_k}\\
\times\frac{(eq^n,fq^n)_k}{(aq^{1+n}/e,aq^{1+n}/f)_k}\left(\frac{q^2a^2}{bcdef}\right)^k\\
=\frac{(aq,aq/ef,\lambda q/e,\lambda q/f)_\infty}
{(aq/e,aq/f,\lambda q/ef,\lambda q)_\infty}\frac{1-\lambda q^{2n}}{1-\lambda}\frac{(\lambda,\lambda b/a,\lambda c/a,\lambda d/a,e,f)_n}{(q,aq/b,aq/c,aq/d,\lambda q/e,\lambda q/f)_n}\left(\frac{aq}{ef}\right)^n\\
\times\sum_{k\geq -n}\frac{(\lambda q^n,\sqrt{\lambda}q^{1+n},-\sqrt{\lambda}q^{1+n},\lambda bq^n/a,\lambda cq^n/a,\lambda dq^n/a)_k}{(q^{1+n},\sqrt{\lambda}q^{n},-\sqrt{\lambda}q^{n},aq^{1+n}/b,aq^{1+n}/c,aq^{1+n}/d)_k}\\
\times\frac{(eq^n,fq^n)_k}{(\lambda q^{1+n}/e,
\lambda q^{1+n}/f)_k}\left(\frac{aq}{ef}\right)^k.
\end{multline}
Next, on both sides of (\ref{a}), replace $a$, $c$, $d$, $e$ and $f$ by $aq^{-2n}$, $cq^{-n}$, $dq^{-n}$, $eq^{-n}$ and $fq^{-n}$ respectively. Note that the condition $\lambda=qa^2/bcd$ implies that $\lambda$ is replaced by $\lambda q^{-2n}$. This yields
\begin{multline*}
\sum_{k\geq -n}\frac{(aq^{-n},q\sqrt{a},-q\sqrt{a},bq^n,c,d,e,f)_k}{(q^{1+n},\sqrt{a},-\sqrt{a},aq^{1-n}/b,aq/c,aq/d,aq/e,aq/f)_k}\left(\frac{q^2a^2}{bcdef}\right)^k\\
=\frac{1-aq^{-2n}}{1-a}\frac{1-\lambda}{1-\lambda q^{-2n}}\frac{(\lambda q^{-2n},\lambda b/a,\lambda cq^{-n}/a,\lambda dq^{-n}/a)_n}{(aq^{-2n},b,cq^{-n},dq^{-n})_n}\left(\frac{a}{\lambda}\right)^{n}\\
\times\frac{(aq^{1-2n},aq/ef,\lambda q/e,\lambda q/f)_\infty}
{(\lambda q^{1-2n},aq/e,aq/f,\lambda q/ef)_\infty}\\
\times\sum_{k\geq -n}\frac{(\lambda q^{-n},q\sqrt{\lambda},-q\sqrt{\lambda},\lambda b q^n/a,\lambda c/a,\lambda d/a,e,f)_k}{(q^{1+n},\sqrt{\lambda},-\sqrt{\lambda},aq^{1-n}/b,aq/c,aq/d,\lambda q/e,
\lambda q/f)_k}\left(\frac{aq}{ef}\right)^k,
\end{multline*}
which is (\ref{semifinite}) after using the simplifications (\ref{elementary1}) and (\ref{elementary3}) on the right-hand side.
\end{proof}
By letting $n\to\infty$ in (\ref{semifinite}), assuming $|qa^2/cdef|<1$ while appealing to Tannery's theorem \cite{Br} for beeing able to interchange limit and summation, one gets the following transformation formula, which was derived in \cite{JS}
\begin{multline}\label{b2}
{}_6\psi_6\!\left[\begin{matrix}q\sqrt{a},-q\sqrt{a},c,d,e,f\\
\sqrt{a},-\sqrt{a},aq/c,aq/d,aq/e,aq/f\end{matrix};
q,\frac{qa^2}{cdef}\right]\\
=\frac{(aq,q/a,aq/ef,aq/cd,\lambda q/e,\lambda
 q/f,aq/\lambda c,aq/\lambda d)_\infty}{(aq/e,aq/f,q/c,q/d,\lambda
 q,q/\lambda,\lambda q/ef,b)_\infty}\\\times
{}_6\psi_6\!\left[\begin{matrix}q\sqrt{\lambda},-q\sqrt{\lambda},
\lambda c/a,\lambda d/a,e,f\\
\sqrt{\lambda},-\sqrt{\lambda},aq/c,aq/d,\lambda
q/e,\lambda q/f\end{matrix};q,\frac{qa^2}{cdef}\right],
\end{multline}
where $\lambda=qa^2/bcd$, and $b$ is now an extra parameter on the right-hand side.\\
As explained in \cite{JS}, an iteration of (\ref{b2}) and an appropriate specialization of both extra parameters appearing on the right-hand side immediately establishes Bailey's formula (\ref{6psi6}).\\

Now, we consider the next level in the hierarchy of identities for very-well-poised nonterminating basic hypergeometric series, which is Bailey's four-term ${}_{10}\phi_9$ transformation~\cite[Appendix (III.39)]{GR} 
\begin{multline}\label{10phi9}
{}_{10}\phi_9\!\left[\begin{matrix}a,\,q\sqrt{a},-q\sqrt{a},b,c,d,e,f,g,h\\
\sqrt{a},-\sqrt{a},aq/b,aq/c,aq/d,aq/e,aq/f,aq/g,aq/h\end{matrix};q,q\right]\\
+\frac{(aq,b/a,c,d,e,f,g,h,bq/c,bq/d,bq/e,bq/f,bq/g,bq/h)_\infty}{(b^2q/a,a/b,aq/c,aq/d,aq/e,aq/f,aq/g,aq/h,bc/a,bd/a,be/a,bf/a,bg/a,bh/a)_\infty}\\
\times{}_{10}\phi_9\!\left[\begin{matrix}b^2/a,\,qb/\sqrt{a},-qb/\sqrt{a},b,bc/a,bd/a,be/a,bf/a,bg/a,bh/a\\b/\sqrt{a},-b/\sqrt{a},bq/a,bq/c,bq/d,bq/e,bq/f,bq/g,bq/h\end{matrix};q,q\right]\\
=\frac{(aq,b/a,\lambda q/f,\lambda q/g,\lambda q/h,bf/\lambda,bg/\lambda,bh/\lambda)_\infty}{(\lambda q,b/\lambda,aq/f,aq/g,aq/h,bf/a,bg/a,bh/a)_\infty}\\
\times{}_{10}\phi_9\!\left[\begin{matrix}\lambda,\,q\sqrt{\lambda},-q\sqrt{\lambda},b,\lambda c/a,\lambda d/a,\lambda e/a,f,g,h\\\sqrt{\lambda},-\sqrt{\lambda},\lambda q/b,aq/c,aq/d,aq/e,\lambda q/f,\lambda q/g,\lambda q/h\end{matrix};q,q\right]\\
+\frac{(aq,b/a,f,g,h,bq/f,bq/g,bq/h,\lambda c/a,\lambda d/a,\lambda e/a,abq/\lambda c,abq/\lambda d,abq/\lambda e)_\infty}{(b^2q/\lambda,\lambda/b,aq/c,aq/d,aq/e,aq/f,aq/g,aq/h,bc/a,bd/a,be/a,bf/a,bg/a,bh/a)_\infty}\\
\times{}_{10}\phi_9\!\left[\begin{matrix}b^2/\lambda,\,qb/\sqrt{\lambda},-qb/\sqrt{\lambda},b,bc/a,bd/a,be/a,bf/\lambda,bg/\lambda,bh/\lambda\\b/\sqrt{\lambda},-b/\sqrt{\lambda},bq/\lambda,abq/\lambda c,abq/\lambda d,abq/\lambda e,bq/f,bq/g,bq/h\end{matrix};q,q\right],
\end{multline}
where $\lambda=qa^2/cde$ and $q^2a^3=bcdefgh$.\\
We can deduce from (\ref{10phi9}) the following semi-finite identity.
\begin{Proposition}
\begin{multline}\label{semifinite10phi9}
\sum_{k\geq -n}\frac{(aq^{-n},q\sqrt{a},-q\sqrt{a},b,cq^n,d,e,f,g,h)_k}{(q^{1+n},\sqrt{a},-\sqrt{a},aq/b,aq^{1-n}/c,aq/d,aq/e,aq/f,aq/g,aq/h)_k}q^k\\
+\alpha_n\;{}_{10}\phi_9\!\left[\begin{matrix}b^2/a,\,bq/\sqrt{a},-bq/\sqrt{a},bq^{-n},bcq^{n}/a,bd/a,be/a,bf/a,bg/a,bh/a\\b/\sqrt{a},-b/\sqrt{a},bq^{1+n}/a,bq^{1-n}/c,bq/d,bq/e,bq/f,bq/g,bq/h\end{matrix};q,q\right]\\
=\beta_n\sum_{k\geq -n}\frac{(\lambda q^{-n},q\sqrt{\lambda},-q\sqrt{\lambda},b,\lambda cq^{n}/a,\lambda d/a,\lambda e/a,f,g,h)_k}{(q^{1+n},\sqrt{\lambda},-\sqrt{\lambda},\lambda q/b,aq^{1-n}/c,aq/d,aq/e,\lambda q/f,\lambda q/g,\lambda q/h)_k}q^k\\
+\gamma_n\;{}_{10}\phi_9\!\left[\begin{matrix}b^2/\lambda,\,bq/\sqrt{\lambda},-bq/\sqrt{\lambda},bq^{-n},bcq^{n}/a,bd/a,be/a,bf/\lambda,bg/\lambda,bh/\lambda\\b/\sqrt{\lambda},-b/\sqrt{\lambda},bq^{1+n}/\lambda,abq^{-n}/\lambda c,abq/\lambda d,abq/\lambda e,bq/f,bq/g,bq/h\end{matrix};q,q\right],
\end{multline}
where 
$\lambda=qa^2/cde$, $c=q^2a^3/bdefgh$, and
\begin{multline*}
\alpha_n=-\frac{b}{a}\frac{(q,q/a,c/b)_n}{(q/b,c/a)_n}\frac{(bq^{1+n}/a,cq^n,aq,bq/c)_\infty}{(bcq^{n}/a,b^2q/a,aq/b,aq/c)_\infty}\\
\times\frac{(bq/d,bq/e,bq/f,bq/g,bq/h,d,e,f,g,h)_\infty}{(bd/a,be/a,bf/a,bg/a,bh/a,aq/d,aq/e,aq/f,aq/g,aq/h)_\infty},
\end{multline*}
\begin{multline*}
\beta_n=\frac{(q/a,\lambda c/a,aq/\lambda d,aq/\lambda e,b/a)_n}{(q/\lambda,c,q/d,q/e,b/\lambda)_n}\\
\times\frac{(aq,bf/\lambda,bg/\lambda,bh/\lambda,\lambda q/f,\lambda q/g,\lambda q/h,bq^n/a)_\infty}{(\lambda q,bf/a,bg/a,bh/a,aq/f,aq/g,aq/h,bq^n/\lambda)_\infty},
\end{multline*}
\begin{multline*}
\gamma_n=\frac{(q,q/a,b/a,aq/\lambda d,aq/\lambda e,\lambda c/ab)_n}{(c,c/a,q/b,q/d,q/e,qb/\lambda)_n}\frac{(bq^n/a,aq,f,g,h)_\infty}{(bcq^{n}/a,b^2q/\lambda,aq/f,aq/g,aq/h)_\infty}\\
\times\frac{(\lambda c/a,\lambda d/a,\lambda e/a,bq/f,bq/g,bq/h,abq/\lambda c,abq/\lambda d,abq/\lambda e)_\infty}{(bd/a,be/a,bf/a,bg/a,bh/a,\lambda/b,aq/c,aq/d,aq/e)_\infty}.
\end{multline*}
\end{Proposition}
\begin{proof}
In the first and the third summations of (\ref{10phi9}), shift the index of summation $k$ by $n$, and replace $a$, $c$, $d$, $e$, $f$, $g$ and $h$ by $aq^{-2n}$, $cq^{-n}$, $dq^{-n}$, $eq^{-n}$, $fq^{-n}$, $gq^{-n}$ and $hq^{-n}$ respectively. Note that the condition $\lambda=qa^2/cde$ implies that $\lambda$ is replaced by $\lambda q^{-2n}$, and the condition $c=q^2a^3/bdefgh$ leaves $c$ unchanged. The first ${}_{10}\phi_9$ is then equal to
$$\delta_n\sum_{k\geq -n}\frac{(aq^{-n},q\sqrt{a},-q\sqrt{a},b,cq^n,d,e,f,g,h)_k}{(q^{1+n},\sqrt{a},-\sqrt{a},aq/b,aq^{1-n}/c,aq/d,aq/e,aq/f,aq/g,aq/h)_k}q^k,$$
where 
\begin{multline*}
\delta_n=\frac{(aq^{-2n},bq^{-n},c,dq^{-n},eq^{-n},fq^{-n},gq^{-n},hq^{-n})_n}{(q,aq^{1-n}/b,aq^{1-2n}/c,aq^{1-n}/d,aq^{1-n}/e,aq^{1-n}/f,aq^{1-n}/g,aq^{1-n}/h)_n}\\
\times\frac{1-a}{1-aq^{-2n}}q^n,
\end{multline*}
and (\ref{10phi9}) is then equivalent to
\begin{multline*}
\sum_{k\geq -n}\frac{(aq^{-n},q\sqrt{a},-q\sqrt{a},b,cq^n,d,e,f,g,h)_k}{(q^{1+n},\sqrt{a},-\sqrt{a},aq/b,aq^{1-n}/c,aq/d,aq/e,aq/f,aq/g,aq/h)_k}q^k\\
+\frac{a_n}{\delta_n}\;{}_{10}\phi_9\!\left[\begin{matrix}b^2/a,\,bq/\sqrt{a},-bq/\sqrt{a},bq^{-n},bcq^{n}/a,bd/a,be/a,bf/a,bg/a,bh/a\\b/\sqrt{a},-b/\sqrt{a},bq^{1+n}/a,bq^{1-n}/c,bq/d,bq/e,bq/f,bq/g,bq/h\end{matrix};q,q\right]\\
=\frac{b_n}{\delta_n}\sum_{k\geq -n}\frac{(\lambda q^{-n},q\sqrt{\lambda},-q\sqrt{\lambda},b,\lambda cq^{n}/a,\lambda d/a,\lambda e/a,f,g,h)_k}{(q^{1+n},\sqrt{\lambda},-\sqrt{\lambda},\lambda q/b,aq^{1-n}/c,aq/d,aq/e,\lambda q/f,\lambda q/g,\lambda q/h)_k}q^k\\
+\frac{c_n}{\delta_n}\;{}_{10}\phi_9\!\left[\begin{matrix}b^2/\lambda,\,bq/\sqrt{\lambda},-bq/\sqrt{\lambda},bq^{-n},bcq^{n}/a,bd/a,be/a,bf/\lambda,bg/\lambda,bh/\lambda\\b/\sqrt{\lambda},-b/\sqrt{\lambda},bq^{1+n}/\lambda,abq^{-n}/\lambda c,abq/\lambda d,abq/\lambda e,bq/f,bq/g,bq/h\end{matrix};q,q\right],
\end{multline*}
where
\begin{multline*}
a_n=\frac{(aq^{1-2n},bq^n/a,c,dq^{-n},eq^{-n},fq^{-n},gq^{-n},hq^{-n})_\infty}{(b^2q/a,aq^{-n}/b,aq^{1-2n}/c,aq^{1-n}/d,aq^{1-n}/e,aq^{1-n}/f,aq^{1-n}/g,aq^{1-n}/h)_\infty}\\
\times\frac{(bq^{1-n}/c,bq/d,bq/e,bq/f,bq/g,bq/h)\infty}{(bcq^n/a,bd/a,be/a,bf/a,bg/a,bh/a)\infty},
\end{multline*}
\begin{multline*}
b_n=\frac{(aq^{1-2n},bq^n/a,\lambda q^{1-n}/f,\lambda q^{1-n}/g,\lambda q^{1-n}/h,bf/\lambda,bg/\lambda,bh/\lambda)_\infty}{(\lambda q^{1-2n},bq^n/\lambda,aq^{1-n}/f,aq^{1-n}/g,aq^{1-n}/h,bf/a,bg/a,bh/a)_\infty}\frac{1-\lambda}{1-\lambda q^{-2n}}q^n\\
\times\frac{(\lambda q^{-2n},bq^{-n},\lambda c/a,\lambda dq^{-n}/a,\lambda eq^{-n}/a,fq^{-n},gq^{-n},hq^{-n})_n}{(q,\lambda q^{1-n}/b,aq^{1-2n}/c,aq^{1-n}/d,aq^{1-n}/e,\lambda q^{1-n}/f,\lambda q^{1-n}/g,\lambda q^{1-n}/h)_n},
\end{multline*}
\begin{multline*}
c_n=\frac{(aq^{1-2n},bq^n/a,fq^{-n},gq^{-n},hq^{-n},bq/f,bq/g,bq/h)_\infty}{(b^2q/\lambda,\lambda q^{-n}/b,aq^{1-2n}/c,aq^{1-n}/d,aq^{1-n}/e,aq^{1-n}/f,aq^{1-n}/g,aq^{1-n}/h)_\infty}\\
\times\frac{(\lambda c/a,\lambda dq^{-n}/a,\lambda eq^{-n}/a,abq^{1-n}/\lambda c,abq/\lambda d,abq/\lambda e)\infty}{(bcq^n/a,bd/a,be/a,bf/a,bg/a,bh/a)\infty}.
\end{multline*}
Using the simplifications (\ref{elementary1}) and (\ref{elementary3}), we get $a_n/\delta_n=\alpha_n$, $b_n/\delta_n=\beta_n$ and $c_n/\delta_n=\gamma_n$, which yields (\ref{semifinite10phi9}).
\end{proof}
Let $n\to\infty$ in (\ref{semifinite10phi9}), assuming $|q^2a^3/bdefgh|=|c|<1$ and $|\lambda c/a|=|aq/de|<1$ while appealing to Tannery's theorem \cite{Br} for beeing able to interchange limit and summation. One gets the following transformation formula
\begin{multline}\label{8psi8}
{}_8\psi_8\!\left[\begin{matrix}q\sqrt{a},-q\sqrt{a},b,d,e,f,g,h\\
\sqrt{a},-\sqrt{a},aq/b,aq/d,aq/e,aq/f,aq/g,aq/h\end{matrix};
q,c\right]\\
=\frac{(aq,q/a,\lambda c/a,aq/\lambda d,aq/\lambda e,b/a,bf/\lambda,bg/\lambda,bh/\lambda,\lambda q/f,\lambda q/g,\lambda q/h)_\infty}{(\lambda q,q/\lambda,c,q/d,q/e,b/\lambda,bf/a,bg/a,bh/a,aq/f,aq/g,aq/h)_\infty}\\\times
{}_8\psi_8\!\left[\begin{matrix}q\sqrt{\lambda},-q\sqrt{\lambda},b,
\lambda d/a,\lambda e/a,f,g,h\\
\sqrt{\lambda},-\sqrt{\lambda},\lambda q/b,aq/d,aq/e,\lambda
q/f,\lambda q/g,\lambda q/h\end{matrix};q,\frac{\lambda c}{a}\right]\\
+\frac{b}{a}\frac{(q,q/a,c/b,aq,bq/c,bq/d,bq/e,bq/f,bq/g,bq/h,d,e,f)_\infty}{(q/b,c/a,b^2q/a,aq/b,aq/c,bd/a,be/a,bf/a,bg/a,bh/a,aq/d,aq/e,aq/f)_\infty}\\
\times\frac{(g,h)_\infty}{(aq/g,aq/h)_\infty}\;{}_{8}\phi_7\!\left[\begin{matrix}b^2/a,\,bq/\sqrt{a},-bq/\sqrt{a},bd/a,be/a,bf/a,bg/a,bh/a\\b/\sqrt{a},-b/\sqrt{a},bq/d,bq/e,bq/f,bq/g,bq/h\end{matrix};q,c\right]\\
+\frac{(q,q/a,b/a,aq/\lambda d,aq/\lambda e,\lambda c/ab,aq,f,g,h)_\infty}{(c,c/a,q/b,q/d,q/e,qb/\lambda,b^2q/\lambda,aq/f,aq/g,aq/h)_\infty}\\
\times\frac{(\lambda c/a,\lambda d/a,\lambda e/a,bq/f,bq/g,bq/h,abq/\lambda c,abq/\lambda d,abq/\lambda e)_\infty}{(bd/a,be/a,bf/a,bg/a,bh/a,\lambda/b,aq/c,aq/d,aq/e)_\infty}\\
\times{}_{8}\phi_7\!\left[\begin{matrix}b^2/\lambda,\,bq/\sqrt{\lambda},-bq/\sqrt{\lambda},bd/a,be/a,bf/\lambda,bg/\lambda,bh/\lambda\\b/\sqrt{\lambda},-b/\sqrt{\lambda},abq/\lambda d,abq/\lambda e,bq/f,bq/g,bq/h\end{matrix};q,\frac{\lambda c}{a}\right],
\end{multline}
where $c=q^2a^3/bdefgh$, $\lambda=qa^2/cde$, $|c|<1$ and $|\lambda c/a|<1$.\\
Note that when $\lambda=a$ or when $b=a$, this identity is trivial.\\

\noindent{\bf Aknowledgments.} We thank Michael Schlosser for his valuable comments and pointing out a mistake in a previous version of this paper.


\begin{thebibliography}{99}

\bibitem{And} G.\ E.\ Andrews,
{\em Applications of basic hypergeometric functions},
SIAM Rev.\ {\bf 16} (1974), 441--484.

\bibitem{AII} R.\ Askey,
{\em The very well poised ${}_6\psi_6$. II},
Proc.\ Amer.\ Math.\ Soc.\ {\bf 90} (1984), 575--579.

\bibitem{AI} R.\ Askey and M.\ E.\ H.\ Ismail,
{\em The very well poised ${}_6\psi_6$},
Proc.\ Amer.\ Math.\ Soc.\ {\bf 77} (1979), 218--222.

\bibitem{Ba1} W.\ N.\ Bailey,
{\em Series of hypergeometric type which are infinite in both directions},
Quart.\ J.\ Math.\ (Oxford) {\bf 7} (1936), 105--115.

\bibitem{Ba2} W.\ N.\ Bailey,
{\em On the basic bilateral hypergeometric series ${}_2\psi_2$},
Quart.\ J.\ Math.\ (Oxford) (2) {\bf 1} (1950), 194--198.

\bibitem{Br} T.\ J.\ l'A.\ Bromwich,
{\em An introduction to the theory of infinite series},
2nd ed., Macmillan, London, 1949.

\bibitem{Ca} A.-L.\ Cauchy,
{\em M\'emoire sur les fonctions dont plusieurs valeurs sont li\'ees
entre elles par une \'equation lin\'eaire, et sur diverses transformations
de produits compos\'es d'un nombre ind\'efini de facteurs},
C.\ R.\ Acad.\ Sci.\ Paris {\bf 17} (1843), 523;
reprinted in {\em Oeuvres de Cauchy}, Ser.~1 {\bf 8},
Gauthier-Villars, Paris (1893), 42--50.

\bibitem{CF} W. \ Y. \ C. Chen and A. \ M. \ Fu,
{\em Semi-Finite Forms of Bilateral Basic Hypergeometric Series}, Proc. Amer. Math. Soc, to appear.

\bibitem{GR} G.\ Gasper and M.\ Rahman,
{\em Basic Hypergeometric Series}, Encyclopedia of Mathematics
And Its Applications 35, Cambridge University Press, Cambridge, 1990.

\bibitem{Jac} C.\ G.\ J.\ Jacobi,
{\em Fundamenta Nova Theoriae Functionum Ellipticarum},
Regiomonti.\ Sump\-tibus fratrum Borntr\"ager,
1829; reprinted in Jacobi's Gesammelte Werke,
vol.~1, (Reimer, Berlin, 1881--1891), pp.~49--239;  
reprinted by Chelsea (New York, 1969); 
now distributed by the Amer.\ Math.\ Soc., Providence, RI.

\bibitem{JS} F.\ Jouhet and M.\ Schlosser,
{\em Another proof of Bailey's ${}_6\psi_6$ summation},
Aequationes Mathematicae, to appear.

\bibitem{Sc6psi6} M.\ Schlosser,
{\em A simple proof of Bailey's very-well-poised ${}_6\psi_6$ summation},
Proc.\ Amer.\ Math.\ Soc.\ {\bf 130} (2001), 1113--1123.

\bibitem{Sc} M.\ Schlosser,
{\em Abel--Rothe type generalizations of Jacobi's triple product identity},
in ``Theory and Applications of Special Functions.
A Volume Dedicated to Mizan Rahman"
(M.\ E.\ H.\ Ismail and E.\ Koelink, eds.), Dev.\ Math., to appear.

\bibitem{Sl} L.\ J.\ Slater,
{\em Generalized hypergeometric functions}, Cambridge University Press,
London/New York, 1966.

\bibitem{SL} L.\ J.\ Slater and A.\ Lakin,
{\em Two proofs of the ${}_6\psi_6$ summation theorem},
Proc.\ Edinburgh Math.\ Soc.\ (2) {\bf 9} (1953--57), 116--121.

\end{thebibliography}
\end{document}